\input amstex
\documentstyle{amsppt}
\loadbold



\let\le=\leqslant
\let\ge=\geqslant
\define\eps{\varepsilon}
\define\eqdef{\overset\text{\rom{def}}\to=}
\define\codim{\operatorname{codim}}
\define\ord{\operatorname{ord}}
\define\supp{\operatorname{supp}}
\define\mult{\operatorname{mult}}
\define\Bir{\operatorname{Bir}}
\define\Pic{\operatorname{Pic}}
\define\Aut{\operatorname{Aut}}
\define\Bas{\operatorname{Bas}}
\define\Cl{\operatorname{Cl}}
\topmatter
\date 17/FEB/97\enddate
\address
Steklov Institute of Mathematics\endgraf
Russian Academy of Sciences
\endaddress
\email    grin\@ mi.ras.ru \endemail
\author M. M. Grinenko
\endauthor
\title
Birational automorphisms of a three-dimensional double
quadric with an elementary singularity
\endtitle
\thanks
This work was carried out with the financial support of the Russian
Foundation for Basic Research (grant no. 96-01-00820) and a grant
for Support of Leading Scientific Schools (no.~96-15-96146). \linebreak
English translation:  H.H.McFaden
\endthanks
\abstract
It is proved that the group of birational
automorphisms of a three-dimensional
double quadric with a singular point arising from a double
point on the branch divisor is a semidirect product of the
free group generated by birational involutions of a special
form and the group of regular automorphisms. The proof is
based on the method of `untwisting' maximal
singularities of
linear systems.

Bibliography: 5 titles.
\endabstract
\endtopmatter
\rightheadtext{Birational automorphisms of a double
quadric}

\document

\subhead Preliminary remarks\endsubhead
In this article we study the group of birational auto-\linebreak morphisms of a
double quadric with a singular point
arising from a double point on\linebreak the branch divisor
(cut out by a fourth-degree hypersurface in~$\Bbb
P^4$). The smooth\linebreak case was investigated in ~\cite{1}.

It is difficult to say a priori how a double point affects the group
$\Bir(V)$ of birational automorphisms of a variety ~$V$. For example, if
$V$ is a three-dimensional smooth quartic, then, as is well known
~\cite{2}, $\Bir(V)$ coincides with the group of automorphisms; however,
this is no longer so for a quartic with a double point ~\cite{3}. On the
other hand, double spaces that are smooth or have double points do
not have non-trivial birational automorphisms (in all higher
dimensions)~\cite{4}. In our case (an exact formulation is given below)
it turns out, as in the case of a double space, that a singular point
does not add anything new to the group ~$\Bir(V)$ (though it is itself
non-trivial).

All our arguments fit in the framework of the method of maximal
singularities. The first rigorous exposition of this method is due to
Iskovskikh and Manin ~\cite{2}. Recently the method was essentially
reworked by Pukhlikov ~\cite{5}, and this has made it possible to give
the arguments a more geometric and explicit character and extend the
area of application while avoiding the `test class' technique. Our
article also has the purpose of applying the new technique to the study
of varieties with singularities.

The author would like to express his gratitude to A. V. Pukhlikov for the
statement of the problem and for very useful discussions, and to V. A.
Iskovskikh for his support and his attention to this work.

\head
\boldsymbol\S\,1. Formulation of the main result
\endhead

We shall be in the following situation everywhere below. Let $F$~ be a
quartic in
~$\Bbb P^4$ that is smooth everywhere except at the point $P$, where it
has an ordinary double singularity. Assume that the smooth quadric
$Q\subset\Bbb P^4$ is such that the divisor $F_Q=Q\cap F$ passes through
~$P$, has there an ordinary double singularity, and does not have other
singular points. We consider a double cover morphism
$\pi\:V\to Q$ that is branched over ~$F_Q$. The variety $V$ has an
ordinary double point, which we also denote by ~$P$, and does not
contain other singularities.

We assume in addition that there are precisely 12~lines on~$Q$, each
passing through $P$ and tangent to $F_Q$ at some other point different
from
~$P$, that there are no lines on ~$Q$ intersecting $F_Q$ only at~$P$,
and that $F$ and $Q$ do not have common lines passing through ~$P$.
These conditions are obviously satisfied for the general quartic ~$F$.

It is not hard to see that $\Pic(V)=\Bbb Z[H]$ for a double quadric $V$,
where $H$ is the inverse image of a hyperplane in~$\Bbb P^4$. In
fact, more is true: $\Cl(V)=\Pic(V)$, that is, every Weil divisor is
linearly equivalent to a Cartier divisor. Indeed, suppose that
$\varphi\:\widetilde{\Bbb P}^4\to\Bbb P^4$ is the composition of a
double cover of $\Bbb P^4$ with a branch over the quartic $F$ and
a subsequent blow-up of the singular point, and let $\widetilde V$
and $\widetilde F$ be the proper smooth inverse images of $V$ and $F$ in
~$\widetilde{\Bbb P}^4$. Then $\widetilde F=\varphi^*(4H)-2E$ and
$\widetilde V=\varphi^*(2H)-E$, where $E$ is the exceptional divisor
of the blow-up, and this implies that $\widetilde V$ and~$\widetilde
F$ are ample. Moreover, $\Pic(\widetilde F)=\Bbb Z\oplus\Bbb Z$. Using
the Lefschetz theorem on hyperplane sections twice, we find that
$\Cl(\widetilde V)=\Bbb Z\oplus\Bbb Z$, from which we get that
$\Cl(V)=\Bbb Z$.

Let $B\subset V$ be a smooth curve such that $\pi\big|_B\:B\to l$ is an
isomorphism onto a line $l\subset Q$. Then $l$ is one of the
following
three types:
\roster
\item"A." $l$ is tangent to $Q_F$ at two smooth points;
\item"B." $l$ passes through $P$ and is tangent to $Q_F$ at another
point;
\item"C." $l\subset Q_F$.
\endroster
The corresponding curves $B$ will be called lines of type A, B, or C on
~$V$.

As follows from arguments in \cite{1}, there is on $V$ an exactly
one-dimensional family of lines with each line intersecting only
finitely many others.

\subhead Construction of a birational automorphism associated with lines
of type A or B on ~$\pmb{V}$
\endsubhead
Let $B\subset V$ be such a line. Then there exists a line $B^*\subset V$
conjugate to $B$ with respect to the involution $\delta\:V\to V$
transposing the sheets of the cover $\pi\:V\to Q$. We note that $B$ and
$B^*$ intersect at two points. The linear system $|-K_V-B-B^*|$ gives a
rational map $g\:V\dasharrow\Bbb P^2$ that lifts to a morphism
$\overline g\:\overline V\to\Bbb P^2$, where $\psi\:\overline V\to V$ is
the composition of a blow-up of the line~$B$ and then of the proper
inverse image of the line ~$B^*$. We denote by $E$ and $E^*$ the
corresponding exceptional divisors on~$\overline V$.

We remark that the general fibre $X$ of the morphism $\overline g$ is an
elliptic curve that intersects $E^*$ in a point. The reflection with
respect to this point (in the sense of the group law) gives an
automorphism of~$X$. This determines on some open subset of ~$\overline
V$ a regular automorphism that lowers to~$V$ as a birational
automorphism
$\tau_B$ (for details on this see ~\cite{1}).

\subhead Birational rigidity\endsubhead
To determine the birational type of $V$ the following definition is
useful (see~\cite{5}).

\definition{Definition}
A pair $(W,Y)$, where $W$ is a projective variety that is
non-singular in codimension~$1$ and $Y$ is a (Weil) divisor on $W$, is
called a {\it test pair\/} if the following conditions hold:
\roster
\item"a)" the linear system $|Y|$ is unfixed;
\item"b)" there exists a number $\alpha=\alpha(W,Y)\in\Bbb R_+$ such
that for every $\beta>\alpha$ with $\beta\in\Bbb Q$ the linear system
$$
|m(Y+\beta K_W)|
$$
is empty for any $m\in\Bbb Z_+$ with $m\beta\in\Bbb Z$.
\endroster

The indicated number $\alpha(W,Y)$ is called the {\it canonical
threshold\/} of the pair.
\enddefinition

The following are examples of test pairs important in practice:
\roster
\item"--" $(\Bbb P^n,L)$, where $L$ is a hyperplane and here
$\alpha(\Bbb P^n,L)=1/(n+1)$;
\item"--" the pair $(W,Y)$, where $\varphi:W\to S$ is a fibring into
Fano
varieties and $Y$ is the inverse image of a very ample divisor on
$S$, and here the canonical threshold is equal to zero;
\item"--" finally, $(V,-K_V)$ makes it possible to describe the group of
birational automorphisms of the variety $V$ of interest to us, and in
our case $\alpha(V,-K_V)$ is clearly equal to $1$.
\endroster

Our main result is stated in the following theorem.

\proclaim{Theorem}
The singular double quadric $V$ described above is a
birationally rigid variety, that is, for any test pair
$(W,Y)$ and any birational map $\psi\:V\to W$ there is a
birational automorphism $\chi\in\Bir(V)$ such that
$$
\alpha\bigl(V,(\psi\circ\chi)^{-1}(Y)\bigr)\le\alpha(W,Y).
$$
\endproclaim

\proclaim{Corollary}
$1)$ $V$ is non-rational and is not birationally isomorphic to any
fibring into conics or into Del Pezzo surfaces.

$2)$ the group $\Bir(V)$ of birational automorphisms is included in the
exact triple
$$
1\longrightarrow*\tau_B\longrightarrow\Bir(V)\longrightarrow
\Aut(V)\longrightarrow1,
$$
where $\Aut(V)$ is the group of biregular automorphisms of $V$ and
$*\tau_B$ is the free product of the involutions $\tau_B$ over all lines
$B\subset V$ of type \rom A or~\rom B.
\endproclaim

\head
\boldsymbol\S\,2. Maximal singularities of a double quadric with a double point
\endhead

\subhead The origin of maximal singularities\endsubhead
As we have already observed (see \S1), \linebreak $\Pic(V)=\Cl(V)=\Bbb Z[H]$, where $H=\pi^*(L)$,
$\pi\:V\to Q$ is our double cover, and $L$ is a hyperplane in~$\Bbb
P^4$.

Let $(W,Y)$ be a test pair, $\chi\:V\dasharrow W$ a birational map, and
$|D|=\nomathbreak \chi^{-1}(|Y|)$ the proper inverse image of the
linear system ~$|Y|$. Obviously, ${|D|\subset|nH|}$ for some
$n\ge1$, and $\alpha(V,D)=n$.

We assume that $n>\alpha(W,Y)$. Then let $\varphi_0\:V_0\to V$ be a
blow-up of the singular point $P\in V$, and let
$D^0=\varphi^{-1}(D)=\varphi^*(nH)-\nu_0E_0$, where $\nu_0\ge0$ and
$E_0$ is the exceptional divisor. We lift $\chi$ to $V_0$ and resolve
the singularities of the birational map
$\chi_0=\chi\circ\varphi\:V_0\dasharrow W$:
$$
\gather
U
\\
\lower-4pt\hbox{${\sssize\psi}$}\swarrow\hskip4mm \searrow
\lower-4pt\hbox{${\sssize\psi_1}$}
\\
V_0\ \ \ \overset\chi_0\to\dasharrow\ \ \ W
\endgather
$$

In this diagram $\psi$ and $\psi_1$ are birational morphisms. It is
clear that there exist open subsets $W'\subset W$ and $U'\subset U$ such
that $\psi_1\big|_{U'}\:U'\to W'$ is an isomorphism,
$U'$~is smooth, and $\codim_W(W\setminus W')\ge2$.

Since $n>\alpha(W,Y)$, we have that $|Y+nK_W|=\varnothing$, which
implies the condition
$|\psi_1^{-1}(Y)+\nomathbreak nK_{W'}|=\nomathbreak \varnothing$.
On the other hand,
$\psi_1^{-1}(Y)=\psi^{-1}(D^0)$, and, moreover, there exist divisors
$E_1,E_2,\dots,E_k\subset U'$ such that
$$
\gather
K_{U'}=\psi^*K_{V_0}+\sum_i\delta(E_i)E_i,
\\
\psi^{-1}(D^0)=n\psi^*H-\nu_0\psi^*E_0-\sum_i\nu_{E_i}(D^0)E_i.
\endgather
$$
Thus,
$$
\alignat1
\varnothing&=|\psi_1^{-1}(Y)+nK_{U'}|
\\
&=\big|(n-\nu_0)\psi^{-1}(E_0)+\sum_{i=0}^k
\bigl(n\delta(E_i)+(n-\nu_0)\nu_{E_i}(E_0)-\nu_{E_i}(D^0)\bigr)E_i\big|.
\tag 1
\endalignat
$$
We now need the following result.

\proclaim{Lemma}
In our notation $\nu_0\le n$.
\endproclaim

\demo{Proof}
We have the chain of morphisms
$V_0@>{\varphi_0}>>V@>\pi>>Q$.
Let $S'=Q\cap\overline{T_PQ}$, a quadratic cone in ~$\Bbb P^3$, and
let $S=\pi^{-1}(S')$ and $S^0=\varphi^{-1}(S)$. We take a generator
$l'\subset S'$ of the cone. Then either $l=\pi^{-1}(l')$ is a rational
curve with a double point at $P$, or $l=l_1\cup l_2$ is a pair of lines
on ~$V$ that intersect at ~$P$ (and at another point). Setting
$l^0=\varphi^{-1}(l)$, or, respectively, $l^0_i=\varphi^{-1}(l_i)$, we
get that
$$
D^0\circ l^0=2n-2\nu_0<0
$$
in the first case, and
$$
D^0\circ l^0_i=n-\nu_0<0
$$
in the second. This means that $S^0\subset\Bas|D^0|$, which contradicts
the fact that ~$|D^0|$ is unfixed.
\enddemo

Thus, we conclude from \thetag{1} that for some $i$
$$
\nu_{E_i}(D^0)-(n-\nu_0)\nu_{E_i}(E_0)>n\delta(E_i).
$$
Bearing in mind that $\nu_{E_i}(E_0)\ge0$, we obtain the next result.

\proclaim{Proposition}
If $n>\alpha(W,Y)$, then $\nu_{E_i}(D^0)>n\delta(E_i)$ for some $i>0$,
that is, the linear system $|D^0|$ has a maximal singularity in the
sense of ~\cite{2}.
\endproclaim

\remark{Remark}
The number $n$ actually is the canonical threshold of the linear
system ~$|D^0|$. Indeed, let $n_1>n$ be some number, and let $m=n_1-
n>0$. Then obvious arguments give us that
$$
|D^0+n_1K_{V_0}|=|-m(\varphi_0^*(H)-E_0)+(n-\nu_0)E_0|=\varnothing.
$$
Thus, the proposition above is none other than the assertion that there
exist maximal singularities in the linear system ~$|D^0|$.
\endremark

\subhead Maximal singularities on the variety $\pmb{V_0}$
\endsubhead
In the spirit of the more modern approach described in ~\cite{5} we have
at this point proved precisely the following:
{\it if $\chi\:V\dasharrow W$ is a birational map\rom,
$|D|=\chi^{-1}|Y|\subset|nH|$\rom, and $n>\alpha(W,Y)$\rom, then there
exists a triple $(\widetilde V,\varphi,T)$\rom, where $T$ is a prime
Weil divisor on ~$\widetilde V$\rom, $T\not\subset\widetilde
V_{\text{\rom{sing}}}$, and
$\varphi\:\widetilde V\to V_0$ is a birational morphism\rom, such that}
$$
\nu_T|D^0|>n\delta(T)
\tag 2
$$
(here, as usual, $\nu_T|D^0|\eqdef\min_{C\in|D^0|}\nu_T(C)$ and
$\delta(T)\eqdef\ord_T(\varphi^*\omega_{V_0}\otimes\omega_{\widetilde
V}^*)$, and $\omega_{V_0}$ and $\omega_{\widetilde V}$ are the canonical
sheaves on $V_0$ and ~$\widetilde V$).

\smallskip
Accordingly, suppose that the triple $(\widetilde V,\varphi,T)$ gives a
maximal singularity of the linear system
$|D^0|\subset|n\varphi_0^*(H)-\nu_0E_0|$, with $n$ the canonical
threshold, and let $B_0=\varphi(T)$ ($B_0$~ is none other than the
centre of the valuation determined by the divisor ~$T$). Two
essentially different cases are possible.

\smallskip\noindent
{\it Case\/} 1.
$\supp B_0\not\subset E_0$. Then, as is not hard to see, we can lower
everything to $V$ and employ the arguments of ~\cite{1} almost
word for word. Thus, for $B=\varphi_0(B_0)$ we have only one
possibility: $\dim B=1$, and here $B$ is necessarily a line of type A or~B.

\smallskip\noindent
{\it Case\/} 2.
$\supp B_0\subset E_0$. This is the case of so-called infinitely near
maximal singularities (of the linear system $|D|$ on ~$V$). It is
proved below that there are no such singularities in our case.

\head
\boldsymbol\S\,3. Elimination of infinitely near singularities
\endhead

\proclaim{Proposition}
In the preceding notation let $|D^0|\subset|n\varphi_0^*(H)-\nu_0E_0|$
be an unfixed linear system, and let $\nu_0\le n$. There do not
exist
triples $(\widetilde V,\varphi,T)$ \rom{(\it{see} ~\S\,2)} such that the
inequality ~\thetag{2} holds if $\supp\varphi(T)\subset E_0$.
\endproclaim

\demo{Proof}
Suppose the contrary. Then there is a (finite) chain of birational morphisms
(see~\cite{5})
$$
\minCDarrowwidth{12mm}
\CD
\widetilde
V=V_N@>\varphi_{N,N-1}>>V_{N-1}@>\varphi_{N-1,N-
2}>>\dots@>\varphi_{i+1,i}>>
V_i@>\varphi_{i,i-1}>>\dots@>\varphi_{1,0}>>V_0
\endCD
\tag 3
$$

\pagebreak\noindent
such that $\varphi(T)=\varphi_{N,N-
1}\circ\dots\circ\varphi_{1,0}(T)\subset E_0$, each morphism
$\varphi_{k,k-1}$ is a blow-up with centre $B_{k-1}$ at a point or
irreducible (possibly singular) curve on $V_{k-1}$, and\linebreak
$\nu_T|D^0|>n\delta(T)$. Let
$E_i=\varphi_{i,i-1}^{-1}(B_{i-1})\subset V_i$ and $E_N=T$. We can
assume that $B_i$ is a point for $i<L$ and a curve for $i\ge L$, where
$0\le L\le N$; furthermore, we can assume that $B_i\subset E_i$ for
any $i$.

Let $\nu_i=\mult_{B_{i-1}}|D|^i$, where $|D|^i$ is the proper
inverse image of the linear system $|D^0|$ on $V_i$. It is clear that
$\nu_1\ge\nu_2\ge\dots\ge\nu_N$. Moreover, it is easy to see that
$\nu_1\le2\nu_0$.

Let $\Gamma$ be the directed graph of the singularity: vertices $i$
and $j$ are joined by an arrow (from $i$ to $j$) if $i>j$ and
$B_{i-1}\subset E_j^{i-1}$. We denote by $p(i,j)$ the number of
different paths from a vertex $i$ to a vertex $j$ when $i>j$;
$p(i,i)\eqdef1$.

As is well known, the numbers $\delta_i=\codim B_i-1$ and $\nu_i$ are
connected by the Noether--Fano--Iskovskikh inequality
(\cite{1}, \cite{2},~\cite{5}):
$$
\sum_{i=1}^Np(N,i)\nu_i>\sum_{i=1}^Np(N,i)\delta_{i-1}.
\tag 4
$$
We remark that the case $\dim B_0=0$, $N=1$, that is, the maximal
singularity is at a point and is realized by the very first blow-up, is
impossible, because then $2n<\nu_1\le2\nu_0$, so that $\nu_0>n$, which
contradicts the lemma in ~\S\,2.

Suppose now that $\dim B_0=1$. Then it follows from \thetag{4} that
$\nu_1>n$. The surface $E_0$ is isomorphic to a quadric in $\Bbb P^3$,
so that $\Pic(E_0)\simeq\Bbb Z\oplus\Bbb Z$. The linear
system $|D^0|$ restricted to $E_0$ has type $(\nu_0,\nu_0)$, $\nu_0\le
n$; on the other hand, $\ord_{B_0}D^0|_{E_0}\ge\nu_1>n$, which is not
possible.

The remaining case, that is, when $\dim B_0=0$ and $N>1$, is very
difficult; its impossibility is proved in the next section.
\enddemo

\head
\boldsymbol\S\,4. The difficult case
\endhead

The `difficulty' of this case is explained by the following
circumstances. Suppose, for example, that we need to eliminate a
maximal singularity of some linear system $|M|$ over a point ~$A$. To do
this, from the point of view of the approach described in ~\cite{5}, we
must do two things: first, show with the help of the
Noether--Fano--Iskovskikh inequalities that two general elements
$M_1,M_2\in|M|$ cut out a curve $C=M_1\circ M_2$ having high
multiplicity at ~$A$; second, find a surface ~$S$ (a so-called test
surface) that is smooth at $A$, does not contain components of $C$, and
is such that $C\circ S<\mult_PC$.

In our case it is convenient to take as test surfaces the elements of
the class $|\varphi_{1,0}^*(H)-E_0|$, but the following lemma indicates
the source of the troubles.

\proclaim{Lemma}
Let $Z$ be an irreducible reduced curve and let $A\in Z\cap E_0$ be
some point. The following conditions are equivalent\rom:
\roster
\item"(i)" $Z\subset S$ for any surface $S\in|\varphi_{1,0}^*(H)-E_0|$
passing through ~$A;$
\item"(ii)" $\pi\circ\varphi_{1,0}(Z)$ is a line on the quadric~$Q$ and
passing through ~$P$.
\endroster
\endproclaim

\demo{Proof}
It suffices to carry out simple local computations in a neighbourhood of
~$P$.
\enddemo

On the other hand, if the curve $Z$ satisfies the condition ~ii) of the
preceding lemma and the linear system $D^0$ has a maximal singularity
over the point $A$, then it is\,always\,true\,that\,$Z\subset\Bas|D^0|$.
Indeed,\,if\,$\varphi_{1,0}(Z)$\,is\,a\,`line',\,that\,is,
$H\circ\varphi_{1,0}(Z)=1$, \linebreak then $D^0\circ Z=n-\nu_0\le n<\nu_1$, which
is possible only if $Z\subset D^0$. But if $\varphi_{1,0}(Z)$ is a
`conic', that is, $H\circ\varphi_{1,0}(Z)=2$, and in addition
$Z\not\subset D^0$, then we get a contradiction:
$D^0\circ Z=2n-2\nu_0\ge\nu_1>n$, whence $\nu_0<n/2$ (recall that
$2\nu_0\ge\nu_1>n$). Thus, it is not possible to use the test
surface method directly in these cases.

We remark that the points on $E_0$ through which `lines' or `conics'
pass, lie on some curve $C_{F_Q}\subset E_0$ that is a
projectivization of the tangent cone to the branch divisor (so that
$C_{F_Q}$ is a conic on $E_0\cong\Bbb P^1\times\Bbb P^1$). The next
assertion is almost obvious.

\proclaim{Proposition}
A linear system cannot contain maximal singularities over a point
$B_0\in E_0\setminus C_{F_Q}$.
\endproclaim

\demo{Proof}
The proof is by contradiction. Let $m=\mult_{B_0}(D^0_1\circ D^0_2)$, where
$D^0_1,D^0_2$ in $|D^0|$ are sufficiently general. The inequality
$m>4n^2$ can be deduced from~\thetag{4} (see~\cite{1},
\cite{2},~\cite{5}). If $S\in|\varphi_{1,0}^*(H)-E_0|$ is a smooth
surface (this is a surface of type~K3), $B_0\in S$, and $S$ is in
general
position with $D^0_1\circ D^0_2$ (such a surface exists according
to the preceding lemma), then we get a contradiction:
$$
m\le D^0_1\circ D^0_2\circ S=n^2(H)^3-\nu_0^2(E_0)^3=4n^2-2\nu_0^2<4n^2.
$$
\enddemo

Suppose now that the linear system $|D^0|$ has a maximal
singularity over a point $B_0\in E_0$ and let $l^0\ni B_0$ be a `conic'
or `line' on ~$V_0$ (in the last case there exists a `line'~$l^{*0}$
that is conjugate to it with respect to the double cover automorphism).
Let $\eps=\mult_{\,l^0}(D^0_1\circ D^0_2)$
(respectively, $\eps^*=\mult_{\,l^{*0}}(D^0_1\circ D^0_2)$). We
would like to get lower estimates for the quantity $\eps$ or for
a combination of it with other `harmful' parameters. But before doing
this we shall see what upper estimates we can count on.

\proclaim{Proposition}
\rom{(i)} Let $l^0$ be a `conic'. Then
$$
\eps\le2n^2.
\tag 5
$$
\rom{(ii)} Let $l^0$ be a `line', and let $\mu=\mult_{\,l^0}|D^0|$. Then
$$
\eps-\eps^*\le4n^2-\nu_0^2-2\mu n-(\nu_0-\mu)^2-(\nu_1-
\mu)^2.
\tag 6
$$
\endproclaim

\demo{Proof}
(i) Suppose that $D_1$,$D_2\in|D|\subset|nH|$ and $S\in|H|$ are in
general position (we are considering everything on ~$V$). Then
$$
4n^2=D_1\circ D_2\circ S\ge\eps l\circ S=2\eps.
$$

(ii) Suppose that $|S^0|\in|H^*-E_0|$ is a general element passing
through the point $B_0$; it can be assumed that
$S^0\cap\Bas|D^0|=l^0\cup l^{*0}$ (here and below the raised symbol ~$*$
in the notation for a  divisor means the complete inverse image on the
corresponding variety, in this case on ~$V_0$). A natural step is to
blow up ~$l^0$. Let $\psi\:V'\to V_0$ be such a blow-up, and let $E'$ be
its exceptional divisor.

\proclaim{Lemma}
$E'\cong\Bbb P^1\times\Bbb P^1$.
\endproclaim

\demo{Proof}
It is not hard to see that the normal sheaf ${\pmb N}_{l^0|V_0}$ can
be represented as the extension
$$
0\longrightarrow{\pmb O}_{l^0}(-2)\longrightarrow{\pmb
N}_{l^0|V_0}
\longrightarrow{\pmb O}_{l^0}\longrightarrow0,
$$
so that the exceptional divisor $E'$ is isomorphic either to the surface
$F_2$ or to the quadric ~$F_0$.

We assume that $E'\cong F_2$ and we let $s$ and $f$ be the classes of the
exceptional section and of the fibre of the surface ~$E'$. We
consider a maximally general element $T\in|H^*-E_0|$ containing the
curve $l^0$ (this is a smooth surface of type~K3). Let $T'|_E=s+\alpha
f$, where $T'$ is the proper inverse image. For two such elements
$T_1$ and~$T_2$ it is obvious that
$$
2=T'_1\circ T'_2\circ E'=(s+\alpha f)^2=\alpha,
$$
so that $T'\big|_E=s+2f$.

Further, local arguments give us easily that there is a sheaf
$|M|\subset|H^*-
E_0|$ on $V_0$ such that:
\roster
\item"1)" $l^0,l^{*0}\subset\Bas|M|$;
\item"2)" a general element of $|M|$ has a double point at the unique
point of intersection of $l^0$ and $l^{*0}$ (this means that
$\pi\circ\varphi_0(M)$ is tangent to $F_Q$ at the point where
$\pi\circ\varphi_0(l^0)$ is tangent to~$F_Q$);
\item"3)" two general elements of $|M|$ do not have tangencies
along~$l^0$.
\endroster
Suppose that $M'|_E=s+\alpha f$ for the proper inverse image of this
sheaf. It is not hard to see that $2=M'\circ T'\circ E'=\alpha$. On the
other hand, $|M'|$ has the fibre over the singular point of~$M$ as a
base curve and the system $|s+f|$ on $F_2$ does not contain
irreducible curves. Thus, $s\subset\Bas|M'|$, that is, general elements
of $|M|$ are tangent along~$l^0$. This contradiction proves the
lemma.
\enddemo

It is not hard to see that the intersection of the cycles is equal to
$$
\split
D'_1\circ D'_2\circ\ S'&=(nH^*-\nu_0E_0^*-\mu E')^2(H^*-E_0^*-E')\\
&=4n^2-2\nu_0^2-2\mu n+2\mu\nu_0-2\mu^2,
\endsplit
\tag 7
$$
where $|S'|$ is the proper inverse image of the system $|S^0|$.

On the other hand, let $\overline f$ be the fibre of $E'$ over the point
$B_0$. It is easy to see that $\overline f\subset\Bas|D'|$ and
$m_{\overline f}=\ord_{\overline f}(D'_1\circ D'_2)\ge(\nu_1-\mu)^2$ for
general elements $D'_1$,~$D'_2$ of the proper inverse image $|D'|$ of
the system $|D^0|$. Moreover, for the blow-up $\psi$ the linear
system $|D'|$ acquires another base curve $Z'$, $\supp(Z')\subset E'$,
which has degree $d$ with respect to the fibres of $E'$, and
$\mu^2+d=\eps$ (Lemma~6.1(ii) in~\cite{5}). Thus,
$$
D'_1\circ D'_2=m_{\overline f}\overline f+Z'+\eps^*{l^*}'+Z,
$$
where ${l^*}'$ is the proper inverse image of $l^{*0}$ and $Z$ is some
curve whose support does not contain $\overline f$, ${l^*}'$, or
components of the curve ~$Z'$.

Since $S'\cap\Bas|D'|={l^*}'$,
$Z'\circ S'|_{E'}=Z'\circ(s+f)\ge d$, $S'\circ\overline f=1$, and
$S'\circ{l^*}'=-1$, it follows that
$$
D'_1\circ D'_2\circ S'\ge(\nu_1-\mu)^2+d-\eps^*.
$$
Comparison of the expressions obtained yields (ii).
\enddemo

Our next problem, as follows from the proposition just proved, is to get
the strongest possible lower estimates for $\eps$ or,
respectively, for $\eps-\eps^*$.

As before, it is possible to choose a realization of our maximal
singularity (the
sequence \thetag{3}) in which each subsequent centre of a
blow-up dominates the preceding one, and we can let $B_0,\dots,B_{L-1}$
be points and $B_L,\dots,B_{N-1}$  (irreducible) curves.

The nature of our subsequent actions depends essentially on the moment
at which the
proper inverse image of the curve $l^0$ `jumps from' the next centre of
a blow-up. Let
$k=\max(0\dots L-1:B_k\in l^k,\ B_{k+1}\notin l^{k+1})$
(as usual, superscripts denote the proper inverse image in the
corresponding step of the chain of blow-ups).

\smallskip\noindent
{\it Case\/} $k=0$.
In this case the first blow-up takes away the curve of interest to us
from the next centre ~$B_1$. We introduce the notation that has now
become standard for the method of maximal singularities. Namely, for
general elements $D^0_1,D^0_2\in|D^0|$ let
$$
D^0_1\circ D^0_2=Z_0+\widetilde Z_0+\eps l^0,
$$
where $Z_0$ is a curve with $l^0\not\subset\supp Z_0$ and no component
of $Z_0$ lies in ~$E_0$; further, $\widetilde Z_0\subset E_0$ is a base
curve (possibly empty) of the system~$|D^0|$, $\widetilde d$~is its
degree with respect to the class $(1,1)$ on ~$E_0$, and
$\widetilde m_{0,i}=\mult_{B_{i-1}}\widetilde Z^{i-1}_0$
for $i=1,\dots,L$ ($\widetilde d=\nomathbreak \widetilde
m_{0,i}=\nomathbreak 0$ if $\widetilde Z_0=\varnothing$); suppose also
that $Z_i\subset E_i$, $1\le i\le L$, are curves such that
$D^i_1\circ D^i_2=(D^{i-1}_1\circ D^{i-1}_2)^i+Z_i$, and let
$d_i=\deg Z_i$ and $m_{i,j}=\mult_{B_{j-1}}Z^{j-1}_i$ for $L\ge
j>i\ge0$.

The quantities introduced are connected by the following relations
\cite{5}:
$$
\left\{
\aligned
\eps+m_{0,1}+\widetilde m_{0,1}&=\nu_1^2+d_1,
\\
m_{0,2}+\widetilde m_{0,2}+m_{1,2}&=\nu_2^2+d_2,
\\
\hdotsfor2
\\
m_{0,L}+\widetilde m_{0,L}+m_{1,L}+\dots+m_{L-1,L}&=\nu_L^2+d_L.
\endaligned
\right.
\tag 8
$$

We introduce the notation $r_i=p(N,i)$ (the numbers $p(i,j)$ were
defined in ~\S\,3). Let us multiply the rows of the system
~\thetag{8} by $r_1,r_2,\dots,r_L$, respectively, and add them. We get
in a standard way (Theorem~7.1 in~\cite{5}) that
$$
r_1\eps+\sum_{i=1}^Lr_i(m_{0,i}+\widetilde m_{0,i})
\ge\sum_{i=1}^Nr_i\nu_i^2.
\tag 9
$$

In place of \thetag{4} we shall need the refined inequality
(it is derived in exactly the same way as~\thetag{4}, except that ~$E_0$ is also taken
into account)
$$
r_0\nu_0+\sum_{i=1}^Nr_i\nu_i>2n\Sigma_0+n\Sigma_1+r_0n,
\tag 10
$$
where $\Sigma_0=\sum_{i=1}^Lr_i$ and $\Sigma_1=\sum_{i=L+1}^Nr_i$.

\pagebreak
Since $\nu_0\le n$, we can decrease $r_0$ somewhat by setting
$r_0=r_1$; this just makes the inequality ~\thetag{10} stronger.

Then since $2\nu_0\ge\nu_1\ge\dots\ge\nu_N$, we obtain in the usual way
the quadratic inequality
$$
2r_0\nu_0+\sum_{i=1}^Nr_i\nu_i^2>
\frac{(2\Sigma_0+\Sigma_1+r_0)^2n^2}{\Sigma_0+\Sigma_1+\frac12r_0}\,,
$$
from which it follows that
$$
\sum_{i=1}^Nr_i\nu_i^2>4n^2\Sigma_0+2r_0n^2-2r_0\nu_0^2
+\frac{\Sigma_1^2n^2}{\Sigma_0+\Sigma_1+\frac12r_0}\,.
\tag 11
$$

The inequalities \thetag{9} and \thetag{11} together yield
$$
r_1\eps+\sum_{i=1}^Lr_i(m_{0,i}+\widetilde m_{0,i})
\ge4n^2\Sigma_0+2r_0n^2-2r_0\nu_0^2
+\frac{\Sigma_1^2n^2}{\Sigma_0+\Sigma_1+\frac12r_0}\,.
\tag 12
$$

As earlier, let $S^0$ be a general element of the system $|H^*-E_0|$. It
is clear that
$$
m_{0,i}+\widetilde m_{0,i}\le m_{0,1}+\widetilde m_{0,1}
\le D_1^0\circ D_2^0\circ S^0\le4n^2-2\nu_0^2.
$$
Then after division by $r_1=r_0$ we get from the inequality \thetag{12}
that
$$
\eps>2n^2,
$$
which takes care of us if $l^0$ is a `conic'.

In the case when $l^0$ is a `line' we leave the term
$m_{0,1}+\widetilde m_{0,1}$ in \thetag{12} and also divide by~$r_1$:
$$
\eps+m_{0,1}+\widetilde m_{0,1}>6n^2-2\nu_0^2
+\frac{\Sigma_1^2n^2}{(\Sigma_0+\Sigma_1+\frac12r_0)r_1}\,.
\tag 13
$$

As in the proof of the inequality \thetag{6}, let
$\psi\:V'\to V_0$ be a blow-up of the `line'~$l^0$. It is easy to
compute that
$$
D'_1\circ D'_2\circ E'=2\mu^2+2\mu n-2\mu\nu_0.
$$
We carry out a more refined computation of the curves:
$$
D'_1\circ D'_2=\eps^*{l^*}'+Z'+Z_0'+\widetilde Z_0'+C+
\text{(some curves)},
$$
where the curve $C$ is a sum of fibres of $E'$ with multiplicities.
Denoting the zero section and the fibre of ~$E'$ as before by $s$ and
$f$, we can assume that
$$
Z'+C=ds+(\nu_1-\mu)^2f+Mf,
$$
where $M\ge0$ is some number. Next, it is clear that
$\widetilde Z_0'\circ E'=\widetilde m_{0,1}$ and that\linebreak $Z_0'\circ E'=m_{0,1}$,
therefore,
$$
2\mu^2+2\mu n-2\mu\nu_0\ge\eps^*-d-(\nu_1-\mu)^2-M
+\widetilde m_{0,1}+m_{0,1}.
$$
On the other hand,
$$
4n^2-2\nu_0^2-2\mu^2-2\mu n+2\mu\nu_0=D'_1\circ D'_2\circ S'
\ge(\nu_1-\mu)^2+d+\widetilde d-\widetilde m_{0,1}+M,
$$
where $S'$ is a general element of $|H^*-E_0^*-E'|$. The last two
inequalities yield
$$
4n^2-2\nu_0^2\ge\eps^*+\widetilde d+m_{0,1}.
$$
Substituting this in \thetag{13} and taking into account that
$\widetilde m_{0,1}\le\widetilde d$ and $\Sigma_0\ge r_1$, we get the
desired estimate
$$
\eps-\eps^*>2n^2+\frac{\Sigma_1^2n^2}
{(\Sigma_0+\Sigma_1+\frac12r_0)\Sigma_0}\,.
$$

\smallskip\noindent
{\it Case\/} $k>0$.
It is a simple matter to see that the preceding inequalities were valid
because $l^1$ jumped from the centre of the blow-up. We shall try to
do something similar also in the case $k>0$, that is, after making some
refinements, we begin everything as it were from the last step where our
curve
still `sits' on the centre of the blow-up.

We make two preliminary observations about the coefficients $r_i$ in the
Noether--Fano--Iskovskikh inequality ~\thetag{10}. First, starting from
the definition of these numbers, we can write that
$$
r_k=\sum_{j\to k}r_j=\sum\Sb j\to k \\ j\le L \endSb r_j
+\sum\Sb j\to k \\ j>L\endSb r_j.
$$
The inequality \thetag{10} permits us easily to decrease $r_k$; let
$$
r_k=\sum\Sb j\to k \\ j\le L\endSb r_j.
$$
Second, we can let $r_0=r_1=\dots=r_k$ (incidentally, it is not hard to
see that this holds automatically for our choice of realization of the
maximal singularity).

Using the notation for the case $k=0$, we write the {\it two\/} systems
of equalities:
$$
\align
&\left\{
\aligned
\eps+m_{0,1}+\widetilde m_{0,1}&=\nu_1^2+d_1,
\\
\eps+m_{0,2}+m_{1,2}&=\nu_2^2+d_2,
\\
\hdotsfor2
\\
\eps+m_{0,k}+m_{k-1,k}&=\nu_k^2+d_k,
\endaligned
\right.
\tag 14
\endalign
$$

\pagebreak
\noindent
and
$$
\align
\vspace{1\jot}
&\left\{
\aligned
\eps+m_{0,k+1}+m_{k,k+1}&=\nu_{k+1}^2+d_{k+1},
\\
m_{0,k+2}+m_{k,k+2}+m_{k+1,k+2}&=\nu_{k+2}^2+d_{k+2},
\\
\hdotsfor2
\\
m_{0,L}+m_{k,L}+m_{k+1,L}+\dots+m_{L-1,L}&=\nu_L^2+d_L.
\endaligned
\right.
\tag 15
\endalign
$$
We note that $\widetilde m_{0,i}=0$ for $i>1$, since $k>0$.

Next, we multiply each row of the system \thetag{14} by $r_k$, while in
the system ~\thetag{15} we multiply the first row by $r_{k+1}$, the
second by $r_{k+2}$, and so on, after which we add the left-hand and
right-hand sides of both systems. The standard trick with
annihilation of the quantities $d_i$ and $m_{i,j}$ still works, despite
the decrease of the coefficient $r_k$, and, using the quadratic
inequality ~\thetag{11} and taking into account that
$\Sigma_0=kr_k+\sum_{i=k+1}^Lr_i$, we get that
$$
\multline\qquad
\eps\biggl(k+\frac{r_{k+1}}{r_k}\biggr)
+\widetilde
m_{0,1}+\sum_{i=1}^km_{0,i}+\sum_{i=k+1}^L\frac{r_i}{r_k}m_{0,i}
\\
>4n^2\biggl(k+\sum_{i=k+1}^L\frac{r_i}{r_k}\biggr)
+2n^2-
2\nu_0^2+\frac{\Sigma_1^2n^2}{(\Sigma_0+\Sigma_1+\frac12r_0)r_k}\,.
\qquad\endmultline
$$
We now observe that
$$
\sum_{i=k+1}^L\frac{r_i}{r_k}(4n^2-m_{0,i})\ge
(4n^2-m_{0,1})\sum\Sb j\to k \\ j\le L \endSb \frac{r_j}{r_k}=4n^2-
m_{0,1},
$$
and for a general $S^0\in|H^*-E_0|$
$$
4n^2-2\nu_0^2=D_1^0\circ D_2^0\circ S^0\ge m_{0,1}+\widetilde m_{0,1};
$$
we get (using $r_k\ge r_{k+1}$) that
$$
\eps(k+1)+\sum_{i=1}^km_{0,i}>
4n^2k+2n^2+\frac{\Sigma_1^2n^2}{(\Sigma_0+\Sigma_1+\frac12r_0)r_k}\,.
\tag 16
$$

\proclaim{Lemma}
$\sum_{i=1}^km_{0,i}\le4n^2-\eps-\eps^*$ \rom(if $l^0$ is
a `conic', then in place of $\eps^*$ substitute another
~$\eps)$.
\endproclaim

\demo{Proof}
Let $S^0\in|H^*-E_0|$ and let $B_0\in S^0$ be a general element. It is
easy to compute that in the $k$th step
$$
D_1^k\circ D_2^k\circ S^k=4n^2-2\nu_0^2-\sum_{i=1}^k\nu_i^2.
$$
On the other hand,
$$
D_1^k\circ D_2^k=Z_0^k+\widetilde Z_0^k
+\eps l^k+\eps^*l^{*k}+Z_1^k+\dots+Z_k
$$
(if $l^0$ is a `conic', then the corresponding term must be thrown out,
of course).

\pagebreak
It is clear that $S^k=H^*-\sum_{i=0}^kE_i^*$, from which $l^k\circ S^k=-
k$, $l^{*k}\circ S^k=0$, $Z_0^k\circ\nomathbreak S^0\ge m_{0,k}\ge0$,
and $\widetilde Z_0^k=\widetilde d-\widetilde m_{0,1}$ (recall that
$\widetilde d=\deg_{E_0}\widetilde Z_0$).

Further, $Z_i^{i+1}\circ S^{i+1}=d_i-m_{i,i+1}$ for $1\le i<k$ and
hence
$$
Z_i^k\circ S^k=Z_i^{i+1}\circ S^{i+1}-Z_i^k\circ
\biggl(\sum_{j=i+2}^kE_j^*\biggr)=d_i-m_{i,i+1}.
$$
Finally, since $Z_k\circ S^k=d_k$, we get that
$$
4n^2-2\nu_0^2-\sum_{i=1}^k\nu_i^2\ge-k\eps+\widetilde d
-\widetilde m_{0,1}+\sum_{i=1}^{k-1}(d_i-m_{i,i+1})+d_k.
$$
Adding both sides of the system \thetag{14} yields
$$
k\eps+\widetilde m_{0,1}+\sum_{i=1}^km_{0,i}=\sum_{i=1}^k\nu_i^2+
\sum_{i=1}^{k-1}(d_i-m_{i,i+1})+d_k.
$$
Comparing these two expressions, we see that
$$
4n^2-2\nu_0^2\ge\sum_{i=1}^km_{0,i}+\widetilde d.
$$
The assertion of the lemma is now obtained from the estimate for
$\widetilde d$:
$$
\align
2\nu_0^2&=D_1^0\circ D_2^0\circ E_0\ge\eps+\eps^*-
\widetilde d
\quad\text{if $l^0$ is a `line'};
\\
2\nu_0^2&=D_1^0\circ D_2^0\circ E_0\ge2\eps-\widetilde d
\quad\text{if $l^0$ is a `conic'}.
\endalign
$$
\enddemo

Everything is now ready for the needed estimates. Accordingly, let $l^0$
be a conic; the inequality~\thetag{16} together with the statement of
the lemma gives us that
$$
(k-1)\eps>4n^2k-2n^2.
$$
The assumption $k=1$ leads immediately to a contradiction: $0>2n^2$.
But if $k>1$, then it is not hard to see that $\eps>2n^2$.

Suppose now that $l^0$ is a line. Then
$$
k\eps>4n^2k-2n^2+\eps^*+
\frac{\Sigma_1^2n^2}{(\Sigma_0+\Sigma_1+\frac12r_0)r_k}\,.
$$
Since $kr_k\le\Sigma_0$, it follows that
$$
\eps>4n^2-\frac{2n^2}k+\frac{\eps^*}k+
\frac{\Sigma_1^2n^2}{(\Sigma_0+\Sigma_1+\frac12r_0)\Sigma_0}\,.
$$
For $k=1$ we get at once the needed estimate
$\eps-\eps^*>2n^2
+\dfrac{\Sigma_1^2n^2}{(\Sigma_0+\Sigma_1+\frac12r_0)\Sigma_0}$\,.
For $k>1$ we must also take into account the inequality
$\eps+\eps^*\le4n^2$.

Thus, we have proved the following result.
\proclaim{Proposition}
\roster
\item"(i)" $\eps>2n^2$ if $l^0$ is a `conic'\rom;
\item"(ii)"
$\eps-\eps^*>2n^2
+\dfrac{\Sigma_1^2n^2}{(\Sigma_0+\Sigma_1+\frac12r_0)\Sigma_0}$ if
$l^0$ is a `line'.
\endroster
\endproclaim

\proclaim{Corollary 1}
The linear system $|D^0|$ cannot have maximal singularities over a point
$B_0$ through which a `conic' passes.
\endproclaim

\demo{Proof}
The case $k=1$ is impossible, as we have just seen; for other values of
$k$ we get a contradiction of the inequality~\thetag{5}.
\enddemo

\proclaim{Corollary 2}
The linear system $|D^0|$ cannot have maximal singularities over a point
$B_0$ through which a `line' passes.
\endproclaim

\demo{Proof}
We use a device described in ~\cite{3}. The inequality ~\thetag{10}
implies that
$$
\Sigma_1>\frac{2n-\nu_1}{\nu_1-n}\Sigma_0+\frac{n-\nu_0}{\nu_1-n}r_0.
$$
Next, we let $\nu_0+\nu_1=3\theta$ with $n/2<\theta\le n$, and
we consider the auxiliary function
$\Lambda(t)=4n^2-2n\mu-t^2-(t-\mu)^2-(2t-\mu)^2$.

We remark that $\Lambda(t)$ is decreasing for $t>n/2$ and that
$\eps-\eps^*\le\Lambda(\theta)$. Moreover,
$$
\Sigma_1>\frac{2n-2\theta}{2\theta-n}\Sigma_0+\frac{n-\theta}{2\theta-
n}r_0.
$$

\proclaim{Lemma}
$\theta>\frac23n$.
\endproclaim

\demo{Proof}
Suppose not. Then $\Sigma_1>2\Sigma_0+r_0$ by the preceding inequality
for $\Sigma_1$, and hence
$$
\eps-\eps^*>2n^2
+\frac{\Sigma_1^2n^2}{(\Sigma_0+\Sigma_1+\frac12r_0)\Sigma_0}
\ge 2n^2+\frac23n^2\frac{\Sigma_1}{\Sigma_0}>\frac{10}3n^2.
$$
On the other hand,
$\eps-\eps^*<\Lambda(\frac n2)\le\frac{21}8n^2$, which is a contradiction.
\enddemo

Thus, $\theta>\frac23n$. But then
$$
2n^2<\eps-\eps^*\le\Lambda\Bigl(\frac23n\Bigr)
\le\frac{33}{18}n^2<2n^2.
$$
The corollary is proved.
\enddemo

\head
\boldsymbol\S\,5. Untwisting automorphisms
\endhead

By now we have established the following: if $\chi\:V\dasharrow W$ is a
birational map onto a test pair $(W,Y)$,
$|D|=\chi^{-1}|Y|\subset|nH|$, and $n>\alpha(W,Y)$, then there is a line
$B\subset V$ of type A or B at which the linear system $|D|$ has a
maximal singularity.

On the other hand, a birational automorphism $\tau_B$ associated in a
natural way with the line ~$B$ was constructed in \S\,1. If we now show
that $(\chi\circ\tau_B)^{-1}|Y|\subset|mH|$ and $m<n$, then this means
that for any such birational automorphism there exist lines $B_1,\dots
B_k$, not necessarily all distinct, such that
$(\chi\circ\tau_{B_1}\circ\dots\circ\tau_{B_k})^{-
1}|Y|$ is contained in $|\alpha(W,Y)H|$.

The next proposition shows that this is indeed the case.

\proclaim{Proposition}
Let $|D|\subset|nH|$ be an unfixed linear system with $n>1$ having a
maximal singularity at a line $B\subset V$ of type ~\rom A or ~\rom B.
Then
$$
\tau_B^{-1}|D|\subset|mH|
$$
for $m<n$.
\endproclaim

\demo{Proof}
Let $B^*$ be the line conjugate to $B$ with respect to the
double cover morphism $\pi:V\to Q$. We set $\nu=\mult_B|D|$ ($\nu>n$ by
assumption) and $\nu^*=\mult_{B^*}|D|$.

Suppose that $B$ is a line of type A. Then the proof of the analogous
assertion
in ~\cite{1} carries over word for word to our situation. Furthermore,
$$
\tau_B^{-1}|D|\subset|(9n-8\nu)H|
$$
and it is clear that $9n-8\nu<n$.

Suppose now that $B$ is a line of type B, $\varphi_0\:V_0\to V$ is a
blow-up of the singular point $P\in V$, $E_0$ is the exceptional
divisor, $B^0=\varphi_0^{-1}(B)$, and $B^{*0}=\varphi_0^{-1}(B^*)$;
further, let $\psi\:\overline V\to V_0$ be a successive blow-up first
of ~$B^0$, and then of the proper inverse image ~$B^{*0}$, and let
$E=\psi^{-1}(B^0)$ and $E^*=\psi^{-1}(B^{*0})$.

We consider a section $Q_H$ of $Q$ by a general hyperplane in $\Bbb
P^4$ that passes through $\pi(B)$, $H=\pi^*(Q_H)$,
$$
H^0=\varphi_0^{-1}(H)=H^*-E_0,
$$
and
$$
\overline H=\psi^{-1}(H^0)=H^*-E_0^*-E-E^*.
$$

The linear system $|\overline H|$ gives a morphism
$\varphi_{|\overline H|}\:\overline V\to\Bbb P^2$ that realizes
$\overline V$ as a fibring into curves of arithmetic genus one. The
birational involution $\overline\tau_B$ is defined, if the fibre is
irreducible, to be a fibrewise map with respect to ~$E^*$ (in the sense
of the group law). Let $Q_P=Q\cap\overline{T_PQ}$. Despite the fact that
a whole open subset of $\overline E_0\cup\overline Q_P$, where
$\overline E_0=\psi^{-1}(E_0)$ and
$\overline Q_P=(\pi\circ\varphi_0\circ\psi)^{-1}(Q_P)$, is composed of
reducible fibres, there is an open set
$U=\overline V\setminus\{\text{a finite set of curves}\}$ on which
$\overline\tau_B$ is biregular.

To see this it suffices to consider the restriction of $\overline\tau_B$
to the surface $\overline S$, where
$$
\overline S=(\pi\circ\varphi_0\circ\psi)^{-1}(Q_L)
$$
and $Q_L$ is the
section of $Q$ by a general hyperplane $L\subset\Bbb P^4$, so that
$Q_L\cap Q_P=l\cup l'$, $l'$ is a line passing through ~$P$, and
$Q_L$ does not contain lines of type A or B other than~$l$.

Under such conditions $\overline S$ is a non-singular K3 surface; the
restriction of $\varphi_{|\overline H|}$ to $\overline S$ gives a
morphism $\overline g\:\overline S\to\Bbb P^1$, the curves
$\overline B=E\cap\overline S$ and $\overline B^{@,*}=E^*\cap\overline
S$
are sections of this morphism, and, moreover, $\psi(\overline B)=B^0$
and $\psi(\overline B^{@,*})=B^{*0}$; finally, all the fibres of
~$\overline g$ except one are irreducible. The unique reducible fibre is
$C_1\cup C_2$, where $(\pi\circ\varphi_0\circ\psi)(C_1)=l'$ and
$C_2=\overline E_0\cap\overline S$. We note that $\Pic(\overline S)$ is
generated by the classes $H_S^*=H^*|_{\overline S}$, $\overline B$,
$\overline B^{@,*}$, $C_1$, and ~$C_2$.

Obviously, every birational automorphism is biregular on $\overline S$
and hence $\overline\tau_B\big|_{\overline S}$ extends to an
automorphism of~$\overline S$, and thus $\overline\tau_B$ is biregular
on
~$\overline V$ outside codimen\-sion two.

Our immediate problem is to define the action of $\overline\tau_B$ on
$\Pic(\overline V)$. Since $\overline S$ is invariant with respect to
$\overline\tau_B$, it suffices to carry out the computations on
~$\Pic(\overline S)$.

Let $X$ be the fibre of the morphism
$\overline g\:\overline S\to\Bbb P^1$. We have the following relations
in ~$\Pic(\overline S)$:
$$
\aligned
&X^2=0;
\\
&H_S^{*2}=4;
\\
&\overline B^2=\overline B^{*2}=C_1^2=C_2^2=-2;
\\
&\overline B\circ\overline B^{@,*}=\overline B\circ X
=\overline B\circ C_2=\overline B^{@,*}\circ X=\overline B^{@,*}\circ
C_2=1;
\\
&\overline B\circ C_1=\overline B^{@,*}\circ C_1=X\circ C_1=X\circ
C_2=0;
\\
&C_1\circ C_2=2.
\endaligned
\tag 17
$$

It is clear that
${\overline\tau_B\big|_{\overline S}}^*(X)=X$ and
${\overline\tau_B\big|_{\overline S}}^*(\overline B^{@,*})=\overline
B^{@,*}$.
Moreover, since $C_1\cup C_2$ is invariant with respect to
${\overline\tau_B\big|_{\overline S}}$, $C_1+C_2\sim X$,
$C_1\circ\overline B^{@,*}=0$, and $C_2\circ\overline B^{@,*}=1$, we
have that
${\overline\tau_B\big|_{\overline S}}^*(C_1)=C_1$ and
${\overline\tau_B\big|_{\overline S}}^*(C_2)=C_2$.

This implies the following relations:
$$
\aligned
&{\overline\tau_B\big|_{\overline S}}^*(\overline B)\circ
{\overline\tau_B\big|_{\overline S}}^*(\overline B)=-2;
\\
&{\overline\tau_B\big|_{\overline S}}^*(\overline B)\circ X=1;
\\
&{\overline\tau_B\big|_{\overline S}}^*(\overline B)\circ\overline
B^*=1;
\\
&{\overline\tau_B\big|_{\overline S}}^*(\overline B)\circ C_1=0;
\\
&{\overline\tau_B\big|_{\overline S}}^*(\overline B)\circ C_2=1.
\endaligned
\tag 18
$$
Let
$$
{\overline\tau_B\big|_{\overline S}}^*(\overline B)
=\alpha X+\beta\overline B+\gamma\overline B^{@,*}+\delta C_1+\eps C_2.
$$
Then the relations \thetag{17} and \thetag{18} imply two possible sets
of coefficients:
$$
\alignat4
\beta&=1,\quad& \gamma&=0,\quad& \eps+\alpha&=0,
\quad\text{and}\quad& \eps&=\delta;
\\
\beta&=-1,\quad& \gamma&=2,\quad& \eps+\alpha&=6,
\quad\text{and}\quad& \eps&=\delta.
\endalignat
$$

We are not interested in the first case, because then
$\overline\tau_B\big|_{\overline S}=\operatorname{id}\big|_{\overline
S}$. In the second case, using the fact that $C_1+C_2\sim X$, we have
that
$$
\aligned
\overline\tau_B\big|_{\overline S}^*(\overline B)
&=6H_S^*-7\overline B-4\overline B^{@,*}-6C_2;
\\
\overline\tau_B\big|_{\overline S}^*(H_S^*)
&=7H_S^*-8\overline B-4\overline B^{@,*}-6C_2.
\endaligned
\tag 19
$$

Thus, $\overline\tau_B$ acts on $\Pic(\overline V)$ as follows:
$$
\aligned
\overline\tau_B^*(H^*)&=7H^*-8E-4E^*-6\overline E_0;
\\
\overline\tau_B^*(E)&=6H^*-7E-4E^*-6\overline E_0;
\\
\overline\tau_B^*(E^*)&=E^*;
\\
\overline\tau_B^*(\overline E_0)&=\overline E_0.
\endaligned
\tag 20
$$

Next, let $|\overline D|$ be the proper inverse image of $|D|$ on
$\overline V$. By assumption,
$$
|\overline D|\subset|nH^*-\nu E-\nu^*E^*-\nu_0\overline E_0|.
$$
From the relations \thetag{20},
$$
\overline\tau_B^{-1}|\overline D|\subset|(7n-6\nu)H^*-(8n-7\nu)E-
(4n-4\nu+\nu^*)E^*-\nu_0\overline E_0|.
$$
Finally, lowering all this to $V$, we see that
$\tau_B^{-1}|D|\subset|(7n-6\nu)H|$ and $7n-6\nu<n$. The
proposition is proved.
\enddemo

\head
\boldsymbol\S\,6. Completion of the proof
\endhead

The last thing remaining for us in the proof of the assertions
formulated in~\S\,1 is to see that there are no relations in the group
$B(V)\subset\Bir(V)$ generated by all the birational involutions of the
form
$\tau_B$. Indeed, the arguments of the preceding section give us the
existence of an exact triple
$$
1\longrightarrow B(V)\longrightarrow\Bir(V)\longrightarrow
\Aut(V)\longrightarrow1.
$$
If we now prove that the linear system $|D|$ cannot simultaneously have
two maximal singularities, then we thereby prove the uniqueness of the
process of `untwisting' any birational automorphism (that is, the
uniqueness of the representation
$$
\chi=g\circ\tau_{B_1}\circ\dots\circ\tau_{B_k},
$$
where
$g\in\Aut(V)$), from which it will follow that $B(V)$ is a free
product of birational involutions.

Suppose that the linear system $|D|$ has maximal singularities at the
lines $B_1$ and~$B_2$; we can single out four cases of their mutual
arrangement:
\roster
\item"a)" $B_1\ne B_2,B_1\cap B_2=\varnothing$;
\item"b)" $B_1\ne B_2,B_1\cap B_2=\text{a point}$,
$B_1\ne\delta(B_2)$, where $\delta$ is the double
cover involution;
\item"c)" $B_1\ne B_2,B_1=\delta(B_2)$;
\item"d)" $B_1=B_2$, that is, $B_2$ lies `over' $B_1$ (the centre of an
infinitely near singularity).
\endroster

All this is actually analyzed in \cite{1}: the cases a) and ~d), along
with b) and ~c) under the condition that the intersection point is not
singular, carry over in general without changes to our situation, and
the remaining two cases we consider may differ at most by the
coefficients in the formulae.

\smallskip\noindent
{\it Case\/} b).
Suppose that $B_1\cap B_2=\text{\{a singular point\}}$. We consider a
linear sub-\linebreak system $|S|\subset|2H|$ on $|V|$ such that
$\Bas|S|\cap\Bas|D|=B_1\cup B_2\cup\text{\{isolated points\}}$; the fact
that it is non-empty is ensured by the Riemann--Roch theorem
($|S|$ is actually concerned with separating the fibres of the
cover $\pi\:V\to Q$). We set \linebreak  $\nu_i=\mult_i|D|>n$, $i=1,2$, and let
$\psi\:\overline V\to V$ be the composition of the blow-ups first of the
singular point, and second of the proper inverse images of $B_1$
and~$B_2$.\linebreak
Then for the proper inverse images on $\overline V$ we have that
$\overline D^2\circ\overline S\ge0$. On the other hand, it is not hard
to compute that
$$
\overline D^2\circ\overline S=8n^2-2\nu_0^2-3\nu_1^2-3\nu_2^2-
2(n-\nu_0)(\nu_1+\nu_2)<0
$$
for $0\le\nu_0\le n$. We have obtained a contradiction.

\smallskip\noindent
{\it Case\/} c).
Suppose that the line $\pi(B_1)=\pi(B_2)$ passes through the
singular point on the branch divisor, $S$ is the proper
inverse image of a general hyperplane section of the quadric $Q$ through
this line, and $\psi\:\overline V\to V$ is a composition of blow-ups as
above. A straightforward computation again gives us that
$$
\overline D^2\circ\overline S=4n^2-2\nu_0^2-\nu_1^2-\nu_2^2
-2(n-\nu_0)(\nu_1+\nu_2)-(\nu_1-\nu_2)^2<0
$$
for $0\le\nu_0\le n$, and this contradicts the requirement
$\overline D^2\circ\overline S\ge0$.

The proof of the theorem in \S\,1 is complete.

\enddocument